\documentclass[11pt]{article}
\usepackage[utf8]{inputenc}
\usepackage{authblk}
\usepackage{amsmath}
\usepackage{amssymb}
\usepackage{amsthm}
\usepackage{mathtools}
\usepackage{tikz}
\usepackage{listings}
\usepackage[left=3.5cm,top=3cm,right=3cm,bottom=3cm,]{geometry}
\usepackage{pgfplots}
\usepgfplotslibrary{fillbetween}
\usepackage[
backend=biber,
style=numeric,
sorting=nyt,
maxbibnames=5
]{biblatex}

\newtheorem{thm}{Theorem}[section]
\newtheorem*{thm*}{Theorem}
\newtheorem{defi}[thm]{Definition}
\newtheorem{lem}[thm]{Lemma}

\newtheorem{prop}[thm]{Proposition}
\newtheorem{cor}[thm]{Corollary}

\newtheorem{question}{Question}

\theoremstyle{remark}
\newtheorem{rem}[thm]{Remark}

\DeclareMathOperator{\conv}{conv}

\DeclarePairedDelimiter\abs{\lvert}{\rvert}
\DeclarePairedDelimiter\scalar{\langle}{\rangle}
\DeclarePairedDelimiter{\set}{\{}{\}}

\DeclarePairedDelimiter\floor{\lfloor}{\rfloor}

\DeclareMathOperator{\C}{\mathbb{C}}

\DeclareMathOperator{\N}{\mathbb{N}}

\DeclareMathOperator{\R}{\mathbb{R}}
\DeclareMathOperator{\Z}{\mathbb{Z}}

\DeclareMathOperator{\imag}{\sqrt{-1}}

\DeclareSymbolFont{extraup}{U}{zavm}{m}{n}
\DeclareMathSymbol{\vardiamondsuit}{\mathalpha}{extraup}{87}

\usepackage{xcolor}

\title{Properties of Ehrhart Polynomials whose Roots Lie on the Critical Line}
\author{Max Kölbl\footnote{e-mail: \texttt{max.koelbl@ist.osaka-u.ac.jp}}}
\affil{Department of Pure and Applied Mathematics, Graduate School of Information Science and Technology, Osaka University, Suita, Osaka 565-0871, Japan}
\date{July 2022}

\begin{document}

\maketitle

\begin{abstract}
    We study a class of polynomials that has all of its roots on the critical line and shares many properties with Ehrhart polynomials.
    Braun showed in \cite{braun2008norm} that the roots of Ehrhart polynomials are bounded quadratically and Higashitani \cite{higashitani2012counterexamples} provided examples for polytopes whose Ehrhart polynomial roots come close to this bound.
    In the case of polynomials of the aforementioned class, we present an improved bound.
    As a side effect this confirms a special case of a conjecture from \cite{braun2008ehrhart}. 
\end{abstract}

\section{Introduction}


Given a lattice polytope $P$, we can define its Ehrhart polynomial $E_P$.
It turns out that this is a very useful invariant that holds information about a variety of interesting properties about $P$.
In particular the two highest degree coefficients, but also the $h^*$-polynomial, encode useful geometric information.
This makes it only natural to ask if the roots of $E_P$ also hold a geometric meaning, and to some extent they do:
given two lattice polytopes, $P_1$ and $P_2$, their Ehrhart polynomials will satisfy $E_{P_1}E_{P_2} = E_{P_1\times P_2}$, so decompositions of $E_P$ can yield information about the decomposability of $P$, if the components of $E_P$ are themselves Ehrhart polynomials.
Identifying whether or not a polynomial is Ehrhart is in general a difficult task, but there are some necessary properties that can be used as tests to disconfirm the Ehrhart property in some cases, and one of these properties is root distribution.

An early conjecture for the root distribution of Ehrhart polynomials was presented in \cite{beck2005coefficients}, namely that for all roots $\alpha$ of an Ehrhart polynomial of degree $d$, $-d\leq \mathrm{Re}(\alpha)\leq d-1$.
While this conjecture is true for degree $\leq 5$ and if the roots are real, it was disproved in general in \cite{higashitani2012counterexamples} and \cite{ohsugi2012smooth}.
That was done by constructing polytopes whose Ehrhart polynomials indicate a quadratic bound rather than a linear one.
The quadratic growth of the bound was confirmed in \cite{braun2008norm}:
for an Ehrhart polynomial of degree $d$, the roots all lie in a disc of radius $d(d-\frac12)$ around the point $-\frac12$ in the complex plane.

An important subclass of lattice polytopes are the \emph{reflexive polytopes} which rose to popularity after Batyrev noticed their connection to string theory in \cite{batyrev1994dual}.
They have the property that the roots of their Ehrhart polynomials are distributed symmetrically across the \emph{critical line} 
\begin{align*}
    CL=\set{z\in\C\colon \mathrm{Re}(z)=-\frac12}.
\end{align*}
In particular, certain classes of reflexive polytopes have all of their Ehrhart roots on $CL$, e.g. standard reflexive simplices, hypercubes of the form $[-1,1]^d$, cross polytopes, and a certain class of symmetric edge polytopes that contains the cross polytopes (\cite{higashitani2017interlacing}).
Prompted by a conjecture in \cite{golyshev2009canonical}, such polytopes have also been studied for their own merit, for example in \cite{hegedus2019ehrhart}, where they were named \emph{$CL$-polytopes}.

Let now $\mathfrak{P}$ denote the class of polynomials of the form
\begin{equation}
    f(z) = b(z)(z^2+z+c_0)(z^2+z+c_1)\cdots (z^2+z+c_m)\label{eqn_factored}
\end{equation}
where the $c_i$ are real numbers $\geq \frac14$ and, for a non-zero constant $a$, $b(z)=a$ if $\deg f$ is even and $b(z)=a(2z+1)$ otherwise.
One can verify that all the roots lie symmetrically across the real axis on $CL$.
We call such polynomials \emph{$CL$-polynomials}.
They contain the class of Ehrhart polynomials of $CL$-polytopes.
More precisely, they are contained in a subclass of $\mathfrak{P}$ whose polynomials admit a \emph{$h^*$-polynomial} whose coefficients are palindromic and non-negative $(\vardiamondsuit)$.
The aim of this study is to better understand this subclass.


In Section 2, we will review some preliminaries on Ehrhart theory and interlacing theory.
In Section 3, we present a smaller, but still quadratic, bound than the more general one found in \cite{braun2008norm}.
\begin{thm*}[Theorem \ref{thm_bound}]
Let $p_0^d(z) = d!\left(\binom{z+d}{d}+\binom{z}{d}\right)$ and let $f$ be a $CL$-polynomial of degree $d$ that satisfies $(\vardiamondsuit)$.
Then for all roots $-\frac12+\alpha\imag$ of $f$, $\alpha\leq\Tilde{\alpha}_d$ where
\begin{align*}
    \Tilde{\alpha}_d=\max\set{ \beta\colon -\frac12+\beta\imag\text{ is a root of }p_0^d }.
\end{align*}
\end{thm*}
In particular, this proves a conjecture by Braun and Develin from \cite{braun2008ehrhart} in the case of $CL$-polynomials.
In the mentioned paper, they showed that the polynomial $\binom{z+d}{d}+\binom{z}{d}$ (itself a $CL$-polynomial) satisfies
\begin{align*}
    \abs{\alpha} = \frac{d^2}{\pi} + O(1)
\end{align*}
where $\alpha$ is the imaginary part of the root with the largest absolute value.
The conjecture is that the imaginary parts of any \emph{Stanley non-negative polynomial} (\emph{SSN-polynomial for short}), i.e., polynomials with non-negative $h^*$-polynomial coefficients, are bounded by $\abs{\alpha}$.

Also we show that every root permitted by this bound, including the extremal ones, are obtained by an appropriate $CL$-polynomial that satisfies $(\vardiamondsuit)$ (Theorem \ref{thm_connected}).
Further, we present a sufficient condition for $(\vardiamondsuit)$ and sets of equations which characterise $CL$-polynomials satisfying $(\vardiamondsuit)$ in low degrees. \\\\
\textbf{Acknowledgements.} I would like to express my gratitude to my advisor, Akihiro Higashitani, for many fruitful conversations and useful feedback.

\section{Preliminaries}

In this section, we will recall the basics of Ehrhart theory and some techniques from the theory of interlacing polynomials.

\subsection{Ehrhart polynomials}

Let $\R^d$ be a vector space with integer lattice $M:=\Z^d$.
A \emph{lattice polytope} $P$ is a set $\conv(S)\subset\R^d$ where $S$ is a finite subset of $M$ and $\conv$ denotes the \emph{convex hull} of a set, i.e., the smallest convex subset of $\R^d$ that contains $S$.
Assume $P$ has full-dimension.
We can assign a function $E_P(k)$ to $P$ as follows.
\begin{align*}
    E_P(k) = \abs{M\cap kP}
\end{align*}
where $kP$ denotes the dilation of $P$ by the non-negative integer factor $k$.
In other words, $E_P$ counts the lattice points of the integer dilations of $P$.
In \cite{ehrhart1962geometrie}, Ehrhart proved that these functions are polynomials, hence, we call them \emph{Ehrhart polynomials}.

Lattice polytopes famously have a ring theoretic interpretation as \emph{Ehrhart rings}.
These are the normalisations of the monoid algebras generated by the lattice points of the polytope, interpreted as monomials (see e.g. \cite{bruns2009polytopes}).
An Ehrhart polynomial can then be reinterpreted as Hilbert function of an Ehrhart ring.
Hence, it makes sense to consider the Hilbert series (or \emph{Ehrhart series} in the language of polytopes)
\begin{align*}
    ES_P(t) = \sum_{k\in\N} E_P(k) t^k = \frac{h^*(t)}{(1-t)^{d+1}}.
\end{align*}
Of specific interest is the function $h^*$ in the numerator, which is itself a polynomial (a proof is given e.g. in \cite{stanely1986enumerative}).
We call it the \emph{$h^*$-polynomial} of $P$.
It is also often called the \emph{$\delta$-vector} of $P$.
In fact, this representation as a rational function exists for \emph{any} polynomial and going forward, we will consider the $h^*$-polynomial of arbitrary polynomials.

If the degree of $E_P$ is known, it can be fully retrieved from its $h^*$-polynomial by performing a change of basis:
\begin{align*}
    E_P(z) = \sum_{i = 0}^d h^*_i \binom{z+d-i}{d}
\end{align*}
where $h^*_i$ refers to the $i$th coefficient of the $h^*$-polynomial.

Both the Ehrhart polynomial and the $h^*$-polynomial can encode information about their lattice polytope.
For example, the degree of $E_P$ is equal to the dimension of $P$, the leading coefficient of $E_P$ is its (normalised) volume, and the second highest degree is half its surface area.
Other important properties are captured by the $h^*$-polynomial:
most famously that palindromic coefficients imply that $P$ is \emph{reflexive}, i.e., its dual polytope is also a lattice polytope \cite{hibi1992note}.
Further, $h^*_0=1$, $h^*_1=\abs{P\cap M}-d-1$, $h^*_d=\abs{P^\circ\cap M}$ where $P^\circ$ denotes the interior of $P$, and $\frac{h^*(1)}{d!}$ is equal to the leading coefficient of $E_P$.
See \cite{beck2007computing} for a comprehensive introduction to Ehrhart theory.
In general, the coefficients of the $h^*$-polynomials are always non-negative.

Classifiying the set of attainable $h^*$-polynomials is a major unsolved problem and discussed in great detail in \cite{higashitani2022characterisation}.
One result from that line of research that will be useful later is the following.
\begin{thm*}[Hibi's Lower Bound Theorem \cite{hibi1991lower}]
Let $P$ be a lattice polytope of dimension $d$ with $h^*$-polynomial $h^*(t) = \sum_{i=0}^d h^*_i t^i$.
Further, suppose that $h^*_d\neq 0$.
Then the inequalities $h^*_1\leq h^*_i$ holds for every $i$.
\end{thm*}

\subsection{Interlacing}

The theory of interlacing has gained traction after it was used in \cite{marcus2013interlacing} and \cite{marcus2015interlacing} to prove the Kadison-Singer problem as well as the existence of bipartite Ramanujan graphs.
In Ehrhart theory, interlacing has been used to show that the Ehrhart polynomials of symmetric edge polytopes from bipartite graphs of the form $B_{2,n}$ and $B_{3,n}$ lie on the critical line \cite{higashitani2017interlacing}.

We shall start by defining the interlacing property.


\begin{defi}
Let $f$ and $g$ be polynomials of degrees $d$ and $d+1$ respectively.
Further, let $R$ be a totally ordered subset of $\C$.
We say that \emph{$f$ $R$-interlaces $g$} or \emph{$f$ and $g$ interlace on $R$} if all the roots $a_1,\ldots,a_d$ of $f$ and $b_1,\ldots,b_{d+1}$ lie on $R$ and satisfy
\begin{align*}
    b_1 \leq a_1 \leq b_2 \leq a_2 \leq \cdots \leq a_d \leq b_{d+1}
\end{align*}
with respect to the ordering on $R$.
\end{defi}


Next we quickly recall some results we need in the course of this paper, most of which come from the extensive work of Fisk \cite{fisk2006polynomials}.


\begin{prop}[Theorem 2.1.10 in \cite{rodriguez2010distribution}]\label{thm_rodriguez}
Let $f$ and $g$ be $CL$-polynomials with degrees $d$ and $d+1$ respectively.
Let $h^*_f$ and $h^*_g$ be their respective $h^*$-polynomials.
Assume $h^*_f$ and $h^*_g$ also have degrees $d$ and $d+1$ and their roots interlace on the unit circle.
Then $f$ $CL$-interlaces $g$.
\end{prop}



\begin{prop}[Lemma 1.10 in \cite{fisk2006polynomials}]\label{thm_fisk-lincomb}
If $f$, $g$ have positive leading coefficients and there is an $h$ which is $\R$-interlaced by both $f$ and $g$, then for all positive $\alpha$ and $\beta$ the linear combination $\alpha f + \beta g$ has all real roots and $\R$-interlaces $h$.
\end{prop}



\begin{prop}[Lemma 1.26 in \cite{fisk2006polynomials}, ``Leibnitz Rule'']\label{thm_fisk-leibnitz}
Suppose that $f$, $f_1$, $g$, $g_1$ are polynomials with positive leading coefficients, and with all real roots.
Assume that $f$ and $g$ have no common roots.
If $f_1$ $\R$-interlaces $f$ and $g_1$ $\R$-interlaces $g$, then
$f_1g_1$ $\R$-interlaces $fg_1+f_1g$ which in turn $\R$-interlaces $fg$, $fg_1$, and $f_1g$.
In particular, $fg_1 + f_1g$ has all real roots.
\end{prop}



\begin{prop}[Corollary 1.41 in \cite{fisk2006polynomials}]\label{thm_fisk-limits}
Suppose that $f_1, f_2,\ldots$ and $g_1, g_2,\ldots$ are sequences of polynomials with all real roots that converge to polynomials f and g respectively.
If $f_n$ and $g_n$ $\R$-interlace for all positive integers $n$, then $f$ and $g$ $\R$-interlace.
\end{prop}



\begin{rem}
All the statements from \cite{fisk2006polynomials} refer to interlacing on the real line, but can be transported to any line of the form $c_1\R + c_2$ for complex numbers $c_1,c_2$ by performing an appropriate affine transformation.
Further, since roots are invariant to scaling of their polynomial, positive leading coefficients can also always be achieved.
\end{rem}



\section{Possible roots of $CL$-polynomials}

Before we can study the non-negativity of the coefficients of the $h^*$-polynomial of $CL$-polynomials, we need a useful preliminary result.


\begin{prop}\label{thm_palindromic}
If $f(z)$ is a $CL$-polynomial, its $h^*$-polynomial has palindromic coefficients.
\end{prop}

\begin{proof}
We can prove this inductively over the degree.

Let $b(z)$ be the factor of $f$ mentioned in Equation \eqref{eqn_factored}.
Its Ehrhart series is either
\begin{align*}
    \frac{1}{1-t}\quad\text{or}\quad\frac{t+1}{(1-t)^2}.
\end{align*}
The $h^*$-polynomials of these two Ehrhart series are already palindromic.
Next, we show that if $f(z)$ is of degree $d$ and has a palindromic $h^*$-polynomial, multiplying $f(z)$ by $z^2+z+c$ where $\R\ni c\geq \frac14$ yields a new palindromic $h^*$-polynomial.

Assuming the identity
\begin{align*}
    \sum_{k=0}^\infty f(k)\ t^k =\frac{h^*(t)}{(1-t)^d},
\end{align*}
where $h^*(t) = \sum_{i=0}^{\floor{\frac{d}{2}}} h^*_i t^i$ we obtain
\begin{align*}
    &\sum_{k=0}^\infty (k^2+k+c)\ f(k)\ t^k \\
    &= t\cdot\left(t\cdot\left(\frac{h^*(t)}{(1-t)^{d+1}}\right)^\prime\right)^\prime + t\cdot\left(\frac{h^*(t)}{(1-t)^{d+1}}\right)^\prime + c\cdot \frac{h^*(t)}{(1-t)^{d+1}}.
\end{align*}
Using this, we get
\begin{align*}
    \sum_{k=0}^\infty (k^2+k+c)&f(k) t^k = \sum_{k=0}^\infty \sum_{i=0}^{\floor{\frac{d}{2}}} h^*_i (k^2+k+c)\left(\binom{k+d-i}{d}+\binom{k+i}{d}\right) \\
    &= \sum_{i=0}^{\floor{\frac{d}{2}}} h^*_i\left(t\cdot\left(t\cdot\left(\frac{h^*(t)}{(1-t)^{d+1}}\right)^\prime\right)^\prime + t\cdot\left(\frac{h^*(t)}{(1-t)^{d+1}}\right)^\prime + c\cdot \frac{h^*(t)}{(1-t)^{d+1}}\right).
\end{align*}
This simplifies to
\begin{equation}\label{eqn_palindromic}
    \sum_{i=0}^{\floor{\frac{d}{2}}} h^*_i\left(\frac{\alpha (t^{d+2-i} + t^{i}) + \beta (t^{d+1-n} + t^{i+1}) + \gamma (t^{d-i} + t^{i+2})}{(1-t)^{d+3}}\right)
\end{equation}
where $\alpha = i^2 + i + c$, $\beta = 2(di - i^2 + d + 1 - c)$, and $\gamma = d^2 - 2di - i + i^2 + d + c$.
Hence, the $h^*$-polynomial of $(z^2+z+c)f(z)$ is also palindromic.
\end{proof}

In the following subsection we will make heavy use of this symmetry.


\subsection{An upper bound for the roots of $CL$-polynomials}

It is now time to move on to discussing $CL$-polynomials which satisfy $(\vardiamondsuit)$.
This property is not always fulfilled.
In fact, we will see that for most $CL$-polynomials, it is not.
Subsection \ref{sec_hyperplane} will also give us an idea that $(\vardiamondsuit)$ might be difficult to characterise in general.
However, it is still possible to make some statements about sufficient and necessary conditions.
In this subsection we will introduce an upper bound on the roots of $CL$-polynomials that satisfy $(\vardiamondsuit)$.

Before we start, however, we need some preparations, starting with some notation.
For a positive integer $d$, we define $m_j^d(z) := z+d-j$ and $b_i^d(z) := \prod_{j=0}^{d-1} m_{i+j}^d(z)$.
Hence, $\binom{z+d-i}{d} = \frac{b_i^d(z)}{d!}$.
Further, let $p_i^d(z):=b_i^d(z)+b_{d-i}^d(z)$ if $b_i^d(z)\neq b_{d-i}^d$ and $p_i^d(z) = b_i^d$ otherwise.
Since $f$ is palindromic, it can be expressed in terms of the $p_i$.
\begin{align*}
    d!\ f(z) = \sum_{i=0}^{\floor{\frac{d}{2}}} h^*_i p_i(z)
\end{align*}
Here $h^*_i$ refers to the $i$-th coefficient of the $h^*$-polynomial.

The following three lemmas focus on the behaviour of the functions $p_i^d$.
This information can be used to better understand the behaviour of $f$.


\begin{lem}\label{thm_value-location}
Let $f$ be a $CL$-polynomial.
Then for every $z_0\in CL$, $f(z_0)\in \R \imag^d$.
\end{lem}

\begin{proof}
We shall write $z_0=-\frac12+a_0\imag$.
Without loss of generality, we can assume that $a_0$ is non-negative.
The $m_j(z_0)$ have then the form $(d-j-\frac12)+a_0\imag$.
Notice that $m_{2d-j-1}(z_0) = (-d+j+1-\frac12)+a_0\imag = -(d-j-\frac12)+a_0\imag$.
In other words, we can reformulate $m_j(z_0)$ and $m_{2d-j-1}(z_0)$ in exponential form as
\begin{align*}
    m_j(z_0) &= r^d_j(z_0)\exp(\imag(\frac{\pi}{2}-\mu^d_j(z_0))) \ \text{and}\\
    m_{2d-j-1}(z_0) &= r^d_j(z_0)\exp(\imag(\frac{\pi}{2}+\mu^d_j(z_0)))
\end{align*}
with $r^d_j>0$ and $0\leq \mu^d_{j}(z_0)\leq \frac{\pi}{2}$.
This yields an exponential representation of $b^d_i(z_0)$
\begin{align*}
    b^d_i(z_0) = \prod_{j=0}^{d-1} r^d_{i+j}(z_0) \exp\left(\imag\left(\frac{d\pi}{2}-\sum_{j=0}^{d-1}\mu^d_{i+j}(z_0)\right)\right)
\end{align*}
and of $b^d_{d-i}(z_0)$
\begin{align*}
    b^d_{d-i}(z_0) &= \prod_{j=0}^{d-1} r^d_{d-i+j}(z_0) \exp\left(\imag\left(\frac{d\pi}{2}+\sum_{j=0}^{d-1}\mu^d_{d-i+j}(z_0)\right)\right) \\
    &= \prod_{j=0}^{d-1} r^d_{2d-(i-j+d)-1}(z_0) \exp\left(\imag\left(\frac{d\pi}{2}+\sum_{j=0}^{d-1}\mu^d_{2d-(i-j+d)-1}(z_0)\right)\right) \\
    &= \prod_{j=0}^{d-1} r^d_{i-j+d-1}(z_0) \exp\left(\imag\left(\frac{d\pi}{2}+\sum_{j=0}^{d-1}\mu^d_{i-j+d-1}(z_0)\right)\right) \\
    &= \prod_{j=0}^{d-1} r^d_{i+j}(z_0) \exp\left(\imag\left(\frac{d\pi}{2}+\sum_{j=0}^{d-1}\mu^d_{i+j}(z_0)\right)\right).
\end{align*}
This admits decompositions $b^d_i(z_0) = x^d_i(z_0)y^d_i(z_0)$ and $b^d_{d-i} = x^d_i(z_0)\overline{y^d_i(z_0})$ where
\begin{align*}
    x^d_i(z_0) = \prod_{j=0}^{d-1} r^d_{i+j}(z_0)\imag^d\quad\quad\text{and}\quad\quad y^d_i(z_0) = \exp\left(-\sum_{j=0}^{d-1}\mu_{i+j}^d(z_0)\right).
\end{align*}
Thus, $p_i(z_0) = x^d_i(z_0)(y^d_i(z_0) + \overline{y^d_i(z_0)})$, which lies in $\R \imag^d$.
\end{proof}



\begin{lem}\label{thm_limit-behaviour}
Regard $\R \imag^d$ as a totally ordered set with $a \imag^d \preceq b \imag^d$ if and only if $a \leq b$.

For every $p_i^d(z)$ there exists a positive unique real number $a_i^d$ such that $p_i^d(-\frac12+a_i^d\imag) = 0$ and $p_i^d(-\frac12+b\imag)\succ 0$ for every $b>a_i^d$.
Further, $i<k$ implies $a_i^d>a_k^d$.
\end{lem}

\begin{proof}
Fix an integer $0\leq i\leq \floor{\frac{d}{2}}$.
Without loss of generality, we assume that $\mathrm{Im}(z) \geq 0$.
We recall the decomposition $p_i^d(z) = x^d_i(z)(y^d_i(z)+\overline{y^d_i(z)})$ from the proof of Lemma \ref{thm_value-location}.
We can see that $x^d_i(z)\succ 0$ is always satisfied, so we only need to check the behaviour of $y^d_i(z)+\overline{y^d_i(z)}$.

Let $0\leq j < d$.
Then $\mathrm{Re}(m_j^d(z))>0$ and $\arg(m_j^d(z)) = \arctan\left(\frac{\mathrm{Im}(m_j^d(z))}{\mathrm{Re}(m_j^d(z))}\right)$, yielding 
\begin{align*}
\mu^d_j(z) = \frac{\pi}{2} - \arctan\left(\frac{\mathrm{Im}(m_j^d(z))}{\mathrm{Re}(m_j^d(z))}\right).
\end{align*}
As $\mathrm{Im}(z)$ grows, $\mu^d_j(z)$ tends towards $0$.
Hence, $y^d_i(z)$ tends towards $1$.
For $j \geq d$ we get an analogous result by virtue of symmetry.
Then $a_i$ is the number for which $\mathrm{Re}(y^d_i(-\frac12+a_i)) = 0$ and for all $b>a$, $\mathrm{Re}(y^d_i(-\frac12+b)) > 0$.

Let now $0\leq i< k \leq \frac{d}{2}$ be integers and let $z_i = -\frac12 + a_i\imag$.
At this position, we have
\begin{align*}
    \arg(y^d_i(z_i)) = -\sum_{j=0}^{d-1}\mu^d_{i+j}(z_i) = -\frac{\pi}{2}
\end{align*}
and for every $z=-\frac12+b\imag$, $b>a_i$, we have
\begin{align*}
    \arg(y^d_i(z)) = -\sum_{j=0}^{d-1}\mu^d_{i+j}(z) > -\frac{\pi}{2}.
\end{align*}
If $i>0$, it is useful to notice that in the sum of the $\mu^d_{i+j}$, some cancellations occur.
We saw earlier that $\mu^d_j(z) = - \mu^d_{2d-j-1}(z)$ for all $z\in CL$.
Thus, if there exists an integer $0\leq \Tilde{j}\leq d-1$ such that $i+\Tilde{j}>d-1$, the sum $\sum_{j=0}^{d-1}\mu^d_{i+j}(z)$ will feature both $\mu^d_{i+\Tilde{j}}(z)$ and its opposite.
Since $k$ is greater than $i$, its corresponding sum will contain more pairs of $\mu^d_j(z)$ that cancel out.
Hence, $\arg(y^d_k(z_i)) > \arg(y^d_i(z_i))$ and $\arg(y^d_k(z)) > \arg(y^d_i(z)) > -\frac{\pi}{2}$ for every $z=-\frac12+b\imag$, $b>a_i$.
Thus, $a_i>a_k$.
\end{proof}



\begin{lem}\label{thm_p0-interlacing}
For any positive integer $d$, $p_0^d$ $CL$-interlaces $p_0^{d+1}$.
\end{lem}

\begin{proof}
By Proposition \ref{thm_rodriguez}, a polynomial of degree $d$ interlaces a polynomial of degree $d+1$ on $CL$ if their $h^*$-polynomials have degrees $d$ and $d+1$ respectively and they interlace on the unit circle.
Since the $h^*$-polynomials of $p_0^d$ and $p_0^{d+1}$ are $1+t^d$ and $1+t^{d+1}$ respectively, their roots are $\exp\left(\frac{(1+2n)\pi}{d}\imag\right)$ and $\exp\left(\frac{(1+2n)\pi}{d+1}\imag\right)$ respectively, where $n$ ranges from $0$ to $d-1$ (resp. $d$).
These roots interlace on the unit circle and hence $p_0^d$ and $p_0^{d+1}$ interlace on $CL$.
\end{proof}


Finally, we may discuss the bound of the roots.


\begin{thm}\label{thm_bound}
Let $f$ be a $CL$-polynomial of degree $d$ that satisfies $(\vardiamondsuit)$.
Then for all roots $-\frac12+\alpha\imag$ of $f$, $\alpha\leq\Tilde{\alpha}_d$ where
\begin{align*}
    \Tilde{\alpha}_d=\max\set{ \beta\colon -\frac12+\beta\imag\text{ is a root of }p_0^d }.
\end{align*}
\end{thm}

\begin{proof}
We prove that that $\Tilde{\alpha}_d$ is in fact a bound.
First of all we notice that in the representation of $f$ in terms of the palindromic basis elements $p_i$, the coefficient of $p_0$ cannot be $0$.
That is because all the other $p_i$ have $0$ and $-1$ as common roots, meaning that without $p_0$, $f$ would not be a $CL$-polynomial.
Notice further that $\Tilde{\alpha}_d$ is precisely $a_0^d$ from Lemma \ref{thm_limit-behaviour}.
We want to look at $f(z_\alpha)$ where $z_\alpha=-\frac12+\alpha\imag$.
From Lemma \ref{thm_limit-behaviour} it follows that $p_0(z_\alpha)=0$ and $p_i(z_\alpha)\neq 0$ for $i>0$.
Further, for $z=-\frac12+\beta\imag$, $\beta>\alpha$, $p_i(z)\neq 0$ for all $i$.
Hence, $z_\alpha$ is the largest root that can be assumed by $f$ and it will be assumed if and only if $f$ is equal to $p_0$.
\end{proof}

The polynomial $p_0$ is not itself an Ehrhart polynomial of any polytope.
Hence it is natural to ask what $CL$-polynomials that are also Ehrhart polynomials have large extremal roots.

\begin{prop}\label{thm_simplices}
The \emph{standard reflexive simplex} of dimension $d$ is defined as
\begin{align*}
    \Delta_{sr}^d = \conv(e_1,e_2,\ldots,e_d,-\sum_{i=1}^d e_i)
\end{align*}
and has $h^*$-polynomial $\sum_{i=0}^d t^i$.
Let $a_{sr}\in\R_{\geq 0}$ denote the number such that $-\frac12+a_{sr}^d\imag$ is the extremal Ehrhart root of $\Delta_{sr}^d$ in the upper half plane.
Then every $CL$-polytope of dimension $d$ at most $9$ whose extremal Ehrhart root in the upper half plane is $-\frac12+\beta\imag$ satisfies $\beta\leq a_{sr}^d$.
\end{prop}

\begin{proof}
The reasoning is similar to that in the proof of Theorem \ref{thm_bound}.
There are two cases: $d\leq 5$ and $5<d\leq 9$.
In the case $d\leq 5$ it can be shown computationally that $a_1^d< a_{sr}^d< a_0^d$, which implies that
\begin{align*}
    (\underbrace{p_0}_{\prec 0} + \underbrace{\sum_{i=1}^{\floor{\frac{d}{2}}}p_i}_{\succ 0})(-\frac12+a_{st}^d) = 0.
\end{align*}
Let us now consider another $CL$-polytope of dimension $d$ with Ehrhart polynomial $\sum_{i=0}^{\floor{\frac{d}{2}}} h^*_i p_i$.
Because, $h^*_0$ is always equal to $1$, and $h^*$ is palindromic, we know that $h^*_d=1$ and Hibi's Lower Bound Theorem holds.
Since further every reflexive polytope has a unique interior point (see \cite{beck2007computing}), we get $h^*_i\geq 1$ for all $1\leq i\leq \floor{\frac{d}{2}}$.
Let $a\geq a_{sr}^d$, then we obtain
\begin{align*}
    \abs{p_0(-\frac12+a\imag)} \leq
    \Bigg\lvert\left(\sum_{i=1}^{\floor{\frac{d}{2}}} p_i\right)(-\frac12+a\imag)\Bigg\rvert \leq
    \Bigg\lvert\left(\sum_{i=1}^{\floor{\frac{d}{2}}} h^*_i p_i\right)(-\frac12+a\imag)\Bigg\rvert
\end{align*}
where equality holds if and only if $a=a_{sr}^d$ and $h^*_i=1$ for every $i$.

In the case $5<d\leq 9$ it can be shown computationally that $a_2^d< a_{sr}^d<a_1^d< a_0^d$, which implies that
\begin{align*}
    (\underbrace{p_0+p_1}_{\prec 0} + \underbrace{\sum_{i=2}^{\floor{\frac{d}{2}}}p_i}_{\succ 0})(-\frac12+a_{st}^d) = 0.
\end{align*}
In particular, it can be computed that both $p_0\prec 0$ and $p_1\prec 0$ hold.
Assume $h^*_1 = k$ for some positive integer $k$.
Due to Hibi's Lower Bound Theorem, $k\leq h^*_i$ for all $2\leq i\leq\floor{\frac{d}{2}}$.
Let $a\geq a_{sr}^d$, then we obtain
\begin{align*}
    &\abs{(p_0 + m p_1)(-\frac12+a\imag)} \leq
    \abs{(m p_0 + m p_1)(-\frac12+a\imag)} \\ \leq
    &\Bigg\lvert\left(\sum_{i=2}^{\floor{\frac{d}{2}}} m p_i\right)(-\frac12+a\imag)\Bigg\rvert \leq
    \Bigg\lvert\left(\sum_{i=2}^{\floor{\frac{d}{2}}} h^*_i p_i\right)(-\frac12+a\imag)\Bigg\rvert
\end{align*}
where equality holds if and only if $a=a_{sr}^d$ and $h^*_i = 1$ for every $i$.
\end{proof}

For higher degrees, the inequalities $a_2^d < a_{sr}^d < a_1 < a_2$ are no longer true and Hibi's Lower Bound Theorem can no longer guarantee that the $h^*_i$ for $i\geq 3$ are large enough to balance out $h^*_3$.
In particular, in degree $10$, for $2\leq m\leq 14$ the polynomial
\begin{align*}
    f(z) = p_0(z) + p_1(z) + m p_2(z) + p_3(z) + p_4(z) + p_5(z)
\end{align*}
is a $CL$-polynomial whose extremal roots have a larger absolute imaginary part than those of the Ehrhart polynomial of $\Delta_{sr}$.
At the same time it seems doubtful that any of the eligible polynomials are associated to $CL$-polytopes.


\begin{rem}\label{rem_comparison}
The following table compares the roots $-\frac12+\beta\imag$ of the standard reflexive simplices with the bounds from Theorem \ref{thm_bound}, \cite{braun2008ehrhart}, and \cite{braun2008norm}.

\begin{tabular}{c|c|c|c|c}
    $d$ & $\Tilde{\alpha}_d$ & $\beta$ & $\frac{d^2}{\pi}$ & $d(d-\frac12)$ \\
    \hline
    $2$ & $0.866$ & $0.645$ & $1.273$ & $3$ \\
    $3$ & $2.398$ & $1.658$ & $2.865$ & $7.5$ \\
    $4$ & $4.603$ & $3.040$ & $5.093$ & $14$ \\
    $5$ & $7.457$ & $4.761$ & $7.958$ & $22.5$ \\
    $6$ & $10.952$ & $6.811$ & $11.459$ & $33$ \\
    $7$ & $15.085$ & $9.186$ & $15.597$ & $45.5$ \\
    $8$ & $19.857$ & $11.882$ & $20.372$ & $60$ \\
    $9$ & $25.267$ & $14.899$ & $25.783$ & $76.5$ \\
    $10$ & $31.313$ & $18.236$ & $31.831$ & $95$ \\
    $20$ & $126.802$ & $69.147$ & $127.324$ & $390$ \\
    $30$ & $285.956$ & $151.904$ & $286.479$ & $885$ \\
    $100$ & $3182.575$ & $1622.493$ & $3183.099$ & $9950$ \\
    $150$ & $7161.449$ & $3627.845$ & $7161.972$ & $22425$ 
\end{tabular}
\end{rem}

\subsection{Connectedness of the set of possible roots}

Now that we know the extremal roots that can be obtained, it is natural to ask, what happens in between.
Hence we shall now define the set of all attainable roots.


\begin{defi}
For any positive integer $d$, we define
\begin{align*}
    \Omega_d = \set{z\in CL\colon \exists p\in\mathfrak{P}\colon p(z) = 0\text{ and $(\vardiamondsuit)$ is satisfied}}.
\end{align*}
\end{defi}


As it turns out, all the roots between the two extremal roots can be obtained by some $CL$-polynomial.
Or in other words, $\Omega_d$ is connected.
To prove this, we need one more lemma.


\begin{lem}\label{thm_leibnitz}
For any positive integer $d$, $p_0^{d+1}$ is $CL$-interlaced by $p_0^{d+1} + (2z+1)p_0^d$.
\end{lem}

\begin{proof}
We start with the case where $d$ is odd.
From Lemma \ref{thm_p0-interlacing} we know that $p_0^d$ $CL$-interlaces $p_0^{d+1}$.
Further, $2z+1$ trivially $CL$-interlaces $(2z+1)^2$.
Since $p_0^{d+1}(-\frac12)$ is not a zero ($b_0^{d+1}$ and $b_{d+1}^{d+1}$ are both positive numbers), it does not share a root with $(2z+1)^2$ and hence by Proposition \ref{thm_fisk-leibnitz}, $(2z+1)(p_0^{d+1}+(2z+1)p_0^d)$ interlaces $(2z+1)p_0^{d+1}$.
Dividing $2z+1$ from both expressions yields the statement.

In the case where $d$ is even, $p_0^{d+1}$ has a root at $-\frac12$ due to symmetry.
The root has multiplicity $1$, because if it had a higher multiplicity, $p_0^d$ would need to have two of its roots at $-\frac12$ as well due to interlacing, but we already saw that this is not the case.
Hence we define polynomials $g_i(z) = z^2+z+\frac14+\epsilon_i$ where $\epsilon_1>\epsilon_2>\cdots$ is a sequence of positive reals that goes to $0$.
Obviously, $g_i$ is $CL$-interlaced by $2z+1$ and almost all $g_i$ have no common root with $p_0^{d+1}$.
Hence, by Proposition \ref{thm_fisk-leibnitz}, $(2z+1)p_0^{d+1} + p_0^dg_i$ interlaces $(2z+1)p_0^{d+1}$.
Using Proposition \ref{thm_fisk-limits}, we receive that $(2z+1)(p_0^{d+1}+(2z+1)p_0^d)$ interlaces $(2z+1)p_0^{d+1}$ again and dividing by $2z+1$ yields the statement.
\end{proof}



\begin{thm}\label{thm_connected}
For every positive integer $d$, $\Omega_d$ is connected.
\end{thm}

\begin{proof}
In the case $d=1$, $\Omega_1=\set{-\frac12}$ is a singleton and hence connected.

Next, consider the case $d=2$.
Up to a scalar, the $h^*$-polynomials of $CL$-polynomials of degree $2$ are given by $t^2+c t+1$.
This corresponds to $CL$-polynomials whose roots are $-\frac{1}{2}\pm\frac{\sqrt{c^2-4c-12}}{2c+4}$.
For $0\leq c\leq 6$, the roots lie on $CL$ and $c=0$ corresponds to the boundary of $\Omega_2$.
The maximal value for $c$ corresponds to roots at $-\frac12$.
In particular, the imaginary part of the roots depend monotonously on $c$, so every value on $\Omega_2$ can be assumed.

The proof for higher degrees $d+1$ can be built inductively.
We start by noticing that multiplying any $CL$-polynomial $f$ with $2z+1$ preserves $(\vardiamondsuit)$:
Similarly to how we derived Equation \eqref{eqn_palindromic} earlier, this operation will map $\frac{t^n}{(1-d)^{d+1}}$ to
\begin{align*}
    \frac{(2(d - n) + 1)t^{n+1} + (2n+1) t^n}{(1-t)^{d+1}},
\end{align*}
which has non-negative coefficients.
Thus, for every $z_0\in\Omega_d$, we can construct a $CL$-polynomial of degree $d+1$ which vanishes at $z_0$.
Since we may assume that $\Omega_d$ is connected and the extremal roots are assumed by $p_0^{d+1}$, we only need to consider roots from $(-\frac12-\Tilde{\alpha}_{d+1},-\frac12-\Tilde{\alpha}_{d})\cup(-\frac12+\Tilde{\alpha}_{d},-\frac12+\Tilde{\alpha}_{d+1})$.
Following the notation of Lemma \ref{thm_limit-behaviour} (but this time remembering the degree), we shall denote by $p_i^d$ the palindromic basis element $d! [\binom{z+i}{d}+\binom{z+d-i}{d}]$ and by $a_i^d$ the root on $CL$ of $p_i^d$ with the largest imaginary part.

It is enough to show two things.
Firstly we need to show that $a_0^d < a_0^{d+1}$ because that implies that
\begin{align*}
    p_0^{d+1}(z_0) \prec 0\quad\text{ and }\quad (2z+1)p_0^{d}(z_0)\succ 0
\end{align*}
for every $z_0\in(-\frac12+\Tilde{\alpha}_d\imag, -\frac12+\Tilde{\alpha}_{d+1}\imag)
$, which means that for the right choice of $c$, $p_0^{d+1} + (2z+1) c p_0^{d} = 0$.
Secondly, we need to show that positive linear combinations of $p_0^{d+1}$ and $(2z+1) p_0^d$ are $CL$-polynomials.
Both statements follow from Lemmas \ref{thm_p0-interlacing} and \ref{thm_leibnitz}, as well as Proposition \ref{thm_fisk-lincomb}.
\end{proof}



\section{Further conditions for $(\vardiamondsuit)$}

In this section we use Equation \eqref{eqn_palindromic} to characterise a subset of $\mathfrak{P}$ that satisfies $(\vardiamondsuit)$.

\begin{prop}\label{thm_sufficient-CL}
Let $f(z)$ be a $CL$-polynomial of degree $d$.
Without loss of generality we may assume that the $c_i$ are ordered by size.
Then $f$ satisfies $(\vardiamondsuit)$ if the $c_i$ satisfy
\begin{align*}
    c_i \leq \left\{\begin{matrix}2i+2, & d\text{ is odd} \\ 2i+1, & d\text{ is even.}\end{matrix}\right.
\end{align*}
\end{prop}

\begin{proof}
The proof works inductively.
The idea is to take a $CL$-polynomial of degree $d^\prime = d-2$ that satisfies $(\vardiamondsuit)$, and multiply it with $z^2+z+c$ where $c$ is chosen so that it preserves non-negativity of the coefficients of the $h^*$-polynomial.
That is in particular the case when the three factors, $\alpha = n^2 + n + c$, $\beta = 2(n(d^\prime - n) + d^\prime + 1 - c)$, and $\gamma = (d^\prime - n)^2 - n + d^\prime + c$, are non-negative.
Since $c$ is positive, $\alpha$ and $\gamma$ are always non-negative.
For $\beta$, the largest possible choice for $c$ is $d^\prime+1$.

To complete the induction, we only have to look at the cases of $d=1$ and $d=2$.
We start with the former.
If $f$ has degree $1$, it is of the form $z+\frac12$ and has $h^*$-polynomial $t+1$.
Thus $c_0\leq 2$, $c_1\leq 4$, $c_2\leq 6$ etc.
If $f$ has degree $2$, it is of the form $z^2+z+c_0$ and has $h^*$-polynomial $c_0t^2 + 2(1-c_0)t + c_0$.
Thus, $c_0\leq 1$, $c_1\leq 3$, $c_2\leq 5$ etc.
\end{proof}

\subsection{Hyperplane equations for $\mathfrak{P}$}\label{sec_hyperplane}


\begin{prop}\label{thm_sufficient2}
Let $f(z)$ be a $CL$-polynomial of degree $d$.

For a fixed degree $d$, denote Vieta's formulas by
\begin{align*}
    v_l^d = \sum_{\substack{S\subseteq C^d\\ \abs{S}=l}}\prod_{c_i\in S} c_i
\end{align*}
where $C^d$ is the set of variables $c_0,\ldots,c_{\floor{\frac{d}{2}}}$ that appear in Equation \eqref{eqn_factored}.

If $f$ satisfies $(\vardiamondsuit)$, then the following two inequalities are satisfied.
\begin{itemize}
    \item[(i)] $v_{\floor{\frac{d}{2}}}^d\geq 0$
    \item[(ii)] \begin{align*}
    \left\{\begin{matrix}
    -dv_{\frac{d}{2}}^d +\sum_{l=0}^{\frac{d}{2}-1} 2^{\frac{d}{2}-l} v_l^d\geq 0 & \text{if $d$ is even} \\
    -(d-2)v_{\frac{d-1}{2}}^d +\sum_{l=0}^{\frac{d-1}{2}-1} 3\cdot 2^{\frac{d-1}{2}-l} v_l^d\geq 0 & \text{if $d$ is odd}
    \end{matrix}
    \right.
\end{align*}
\end{itemize}
\end{prop}

\begin{proof}
The two inequalities respectively mean that the first / last and the second / second to last coefficient is non-negative.

If $f$ has the form $a_0+a_1z+\ldots+a_dz^d$, its $h^*$-polynomial has coefficients $h^*_0 =a_0$ and $h^*_1 = -da_0 + \sum_{i=1}^d a_1$.
Without loss of generality, we assume that $f$ is monic.
From this, it can be easily seen that $h^*_0=\prod_{i=0}^{\floor{\frac{d}{2}}} c_i = v_{\floor{\frac{d}{2}}}^d$.

For $h^*_1$, we begin with the case where $d$ is even.
Simplifying Equation \eqref{eqn_factored}, we get terms of the form
\begin{align*}
    c_{i_1}\cdots c_{i_l}z^{2m+n}
\end{align*}
where $2(l+m+n)=d$, $l,m,n\geq 0$.
Hence, $v_l^d$ appears $2^{\frac{d}{2}-l}$ times in the coefficients.
In particular, $v_{\frac{d}{2}}^d$ only appears once as $a_0$.

If $d$ is odd, we also have the term $2z+1$ to account for.
Simplifying Equation \eqref{eqn_factored} yields two types of terms.
\begin{align*}
    2c_{i_1}\cdots c_{i_l}z^{2m+n+1}\quad\text{and}\quad c_{i_1}\cdots c_{i_l}z^{2m+n}
\end{align*}
where $2(l+m+n)=d$, $l,m,n\geq 0$.
Hence, $v_l^d$ appears $3\cdot 2^{\frac{d-1}{2}-l}$ times in the coefficients.
In particular, $v_{\frac{d-1}{2}}^d$ appears once as $a_0$ and twice in the form of $a_1 = 2v_{\frac{d-1}{2}}^d+v_{\frac{d-1}{2}-1}^d$.
Hence, $v_{\frac{d-1}{2}}^d$ appears $d-2$ times.
\end{proof}

Since the $c_i$ are all assumed to be greater or equal to $\frac14$, (i) is trivially satisfied.


A core idea of this proof is the fact that the coefficients of $f$ are just linear combinations of Vieta's formulas.
Since the coefficients of the $h^*$-polynomials are in turn linear combinations of the coefficients of $f$, it follows that the set of polynomials of the form of Equation \eqref{eqn_factored} which satisfy $(\vardiamondsuit)$ form a polyhedral set with respect to Vieta's formulas if the $c_i$ are allowed to take on arbitrary complex values.
Each of the vertices of this polyhedral set corresponds to a polynomial $\binom{z+i}{d}+\binom{z+d-i}{d}$, of which there are exactly $\floor{\frac{d}{2}}$ many.
Hence, we get the following.

\begin{cor}\label{thm_hyperplanes}
The set of degree $d$ polynomials of the form of Equation \eqref{eqn_factored} (allowing the $c_i$ assume arbitrary complex values) parametrised by Vieta's formulas as given in Proposition \ref{thm_sufficient2} is a $\floor{\frac{d}{2}}$-dimensional simplex.
\end{cor}

One should note, however, that not every point in these simplices corresponds to a $CL$-polynomial since the $c_i$ are not required to be real and $\geq \frac14$.

Using \texttt{SAGEMATH} \cite{sagemath}, we were able to compute the hyperplanes for degrees $4$ through $14$ explicitly.
They can be found in the Appendix.
We were also able to computationally verify that the simplices of these degrees are lattice simplices.

\begin{figure}
    \centering
    \begin{tikzpicture}
    \begin{axis}[]
    \addplot [name path=line, domain=1/2:22, samples=2, blue] {(x-0.25)/4};
    \addplot [name path=parabola, domain=-5:6.92, samples=32, red] {x^2/4};
    \path [name path=axis] (axis cs:-5,12) -- (axis cs:6.92,12);
    \path [name path=axis2] (axis cs:-5,-21/16) -- (axis cs:22,-21/16);
    \draw [blue] (axis cs:0.5,1/16) -- (axis cs:0.5,12);
    
    \addplot [black,domain=-2:4,samples=2] {0};
    \addplot [black,domain=4:22,samples=2] {2*(x-4)/3};
    \addplot [black,domain=-2:22,samples=2] {(x+2)/2};
    
    \addplot [
        thick,
        color=blue,
        fill=blue, 
        fill opacity=0.05
    ]
    fill between[
        of=line and axis2,
        soft clip={domain=0.5:22},
    ];
    \addplot [
        thick,
        color=blue,
        fill=blue, 
        fill opacity=0.05
    ]
    fill between[
        of=axis and axis2,
        soft clip={domain=-5:0.5},
    ];
    \addplot [
        thick,
        color=red,
        fill=red, 
        fill opacity=0.05
    ]
    fill between[
        of=parabola and axis,
        soft clip={domain=-5:12},
    ];
    
    \end{axis}
    \end{tikzpicture}

    \caption{The triangle of polynomials in the form of Equation \eqref{eqn_factored} of degree $4$ (including those whose $c_i$ assume arbitrary complex values), with the sets of ``forbidden'' points that don't correspond to $CL$-polynomials.
    In the parabolic area, the corresponding $c_0$ and $c_1$ are complex.
    In the other area, either one or both $c_i$ are less than $\frac14$.}
    \label{fig_poly}
\end{figure}
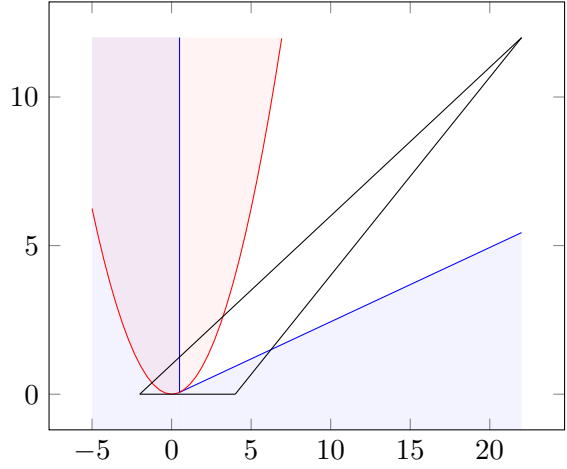

We conclude by posing some questions.

\begin{question}
Are the standard symmetric simplices the $CL$-polytopes whose Ehrhart polynomial roots have the highest attainable imaginary part, even in dimension $\geq 10$?
\end{question}

\begin{question}
How can the simplices from Corollary \ref{thm_hyperplanes} be described for all degrees?
\end{question}

\begin{question}
Are the simplices from Corollary \ref{thm_hyperplanes} always lattice simplices?
\end{question}

\printbibliography

\appendix 

\section{Appendix}

Let $f(z)$ be a $CL$-polynomial of degree $4\leq d\leq 14$ and let Vieta's formulas be defined like in Proposition \ref{thm_sufficient2} for a fixed degree $d$.
Let $\mathbf{v}^d$ denote the vector $\left(v_1^d,\ldots,v_{\floor{\frac{d}{2}}}^d\right)^{tr}$.
Let $\scalar{v,w}$ denote the standard scalar product of $v$ and $w$.

Then $f(z)$ satisfies $(\vardiamondsuit)$ if and only if one of the following holds.

\begin{itemize}
\item The degree of $f$ is 4 and the following inequalities are satisfied.
\begin{itemize}
\item[(i)] $\scalar{(0, 1),\mathbf{v}^{4}}\geq 0$
\item[(ii)] $\scalar{(1, -2),\mathbf{v}^{4}}+2\geq 0$
\item[(iii)] $\scalar{(-2, 3),\mathbf{v}^{4}}+8\geq 0$
\end{itemize}
In particular, $\mathbf{v}^{4}$ lies in the convex hull of

$\set{(22, 12)^{tr}, (4, 0)^{tr}, (-2, 0)^{tr}}$.
\item The degree of $f$ is 5 and the following inequalities are satisfied.
\begin{itemize}
\item[(i)] $\scalar{(0, 1),\mathbf{v}^{5}}\geq 0$
\item[(ii)] $\scalar{(2, -1),\mathbf{v}^{5}}+4\geq 0$
\item[(iii)] $\scalar{(-3, 1),\mathbf{v}^{5}}+54\geq 0$
\end{itemize}
In particular, $\mathbf{v}^{5}$ lies in the convex hull of

$\set{(58, 120)^{tr}, (18, 0)^{tr}, (-2, 0)^{tr}}$.
\item The degree of $f$ is 6 and the following inequalities are satisfied.
\begin{itemize}
\item[(i)] $\scalar{(0, 0, 1),\mathbf{v}^{6}}\geq 0$
\item[(ii)] $\scalar{(2, 1, -3),\mathbf{v}^{6}}+4\geq 0$
\item[(iii)] $\scalar{(8, -8, 15),\mathbf{v}^{6}}+160\geq 0$
\item[(iv)] $\scalar{(-6, 3, -5),\mathbf{v}^{6}}+96\geq 0$
\end{itemize}
In particular, $\mathbf{v}^{6}$ lies in the convex hull of

$\set{(127, 822, 360)^{tr}, (52, 72, 0)^{tr}, (7, -18, 0)^{tr}, (-8, 12, 0)^{tr}}$.
\item The degree of $f$ is 7 and the following inequalities are satisfied.
\begin{itemize}
\item[(i)] $\scalar{(0, 0, 1),\mathbf{v}^{7}}\geq 0$
\item[(ii)] $\scalar{(12, 6, -5),\mathbf{v}^{7}}+24\geq 0$
\item[(iii)] $\scalar{(28, -6, 3),\mathbf{v}^{7}}+296\geq 0$
\item[(iv)] $\scalar{(-96, 12, -5),\mathbf{v}^{7}}+4128\geq 0$
\end{itemize}
In particular, $\mathbf{v}^{7}$ lies in the convex hull of

$\set{(244, 3708, 5040)^{tr}, (118, 600, 0)^{tr}, (34, -72, 0)^{tr}, (-8, 12, 0)^{tr}}$.
\item The degree of $f$ is 8 and the following inequalities are satisfied.
\begin{itemize}
\item[(i)] $\scalar{(0, 0, 0, 1),\mathbf{v}^{8}}\geq 0$
\item[(ii)] $\scalar{(4, 2, 1, -4),\mathbf{v}^{8}}+8\geq 0$
\item[(iii)] $\scalar{(36, 0, -3, 7),\mathbf{v}^{8}}+288\geq 0$
\item[(iv)] $\scalar{(36, -18, 15, -28),\mathbf{v}^{8}}+4824\geq 0$
\item[(v)] $\scalar{(-224, 32, -20, 35),\mathbf{v}^{8}}+9344\geq 0$
\end{itemize}
In particular, $\mathbf{v}^{8}$ lies in the convex hull of

$\set{(428, 13324, 52272, 20160)^{tr}, (232, 3132, 2880, 0)^{tr}, (92, 52, -480, 0)^{tr},\\ (8, -116, 192, 0)^{tr}, (-20, 108, -144, 0)^{tr}}$.
\item The degree of $f$ is 9 and the following inequalities are satisfied.
\begin{itemize}
\item[(i)] $\scalar{(0, 0, 0, 1),\mathbf{v}^{9}}\geq 0$
\item[(ii)] $\scalar{(24, 12, 6, -7),\mathbf{v}^{9}}+48\geq 0$
\item[(iii)] $\scalar{(84, 6, -3, 2),\mathbf{v}^{9}}+600\geq 0$
\item[(iv)] $\scalar{(1188, -126, 27, -14),\mathbf{v}^{9}}+41256\geq 0$
\item[(v)] $\scalar{(-1620, 90, -15, 7),\mathbf{v}^{9}}+137160\geq 0$
\end{itemize}
In particular, $\mathbf{v}^{9}$ lies in the convex hull of

$\set{(700, 39708, 341136, 362880)^{tr}, (412, 11772, 35280, 0)^{tr}, (196, 1404, -3600, 0)^{tr},\\ (52, -468, 720, 0)^{tr}, (-20, 108, -144, 0)^{tr}}$.
\item The degree of $f$ is 10 and the following inequalities are satisfied.
\begin{itemize}
\item[(i)] $\scalar{(0, 0, 0, 0, 1),\mathbf{v}^{10}}\geq 0$
\item[(ii)] $\scalar{(8, 4, 2, 1, -5),\mathbf{v}^{10}}+16\geq 0$
\item[(iii)] $\scalar{(1120, 128, -8, -16, 45),\mathbf{v}^{10}}+7424\geq 0$
\item[(iv)] $\scalar{(920, -26, -4, 7, -15),\mathbf{v}^{10}}+20632\geq 0$
\item[(v)] $\scalar{(272, -224, 68, -56, 105),\mathbf{v}^{10}}+442624\geq 0$
\item[(vi)] $\scalar{(-4520, 260, -50, 35, -63),\mathbf{v}^{10}}+378320\geq 0$
\end{itemize}
In particular, $\mathbf{v}^{10}$ lies in the convex hull of

$\set{(1085, 104008, 1757196, 5132880, 1814400)^{tr}, (680, 36508, 277056, 201600, 0)^{tr},\\ (365, 7528, -3924, -25200, 0)^{tr}, (140, -572, -3024, 7200, 0)^{tr},\\ (5, -392, 2556, -3600, 0)^{tr}, (-40, 508, -2304, 2880, 0)^{tr}}$.
\item The degree of $f$ is 11 and the following inequalities are satisfied.
\begin{itemize}
\item[(i)] $\scalar{(0, 0, 0, 0, 1),\mathbf{v}^{11}}\geq 0$
\item[(ii)] $\scalar{(16, 8, 4, 2, -3),\mathbf{v}^{11}}+32\geq 0$
\item[(iii)] $\scalar{(5904, 792, 36, -42, 35),\mathbf{v}^{11}}+37728\geq 0$
\item[(iv)] $\scalar{(4704, 48, -24, 8, -5),\mathbf{v}^{11}}+85440\geq 0$
\item[(v)] $\scalar{(19216, -1192, 124, -28, 15),\mathbf{v}^{11}}+1740512\geq 0$
\item[(vi)] $\scalar{(-31968, 936, -72, 14, -7),\mathbf{v}^{11}}+4692384\geq 0$
\end{itemize}
In particular, $\mathbf{v}^{11}$ lies in the convex hull of

$\set{(1610, 244708, 7272216, 44339040, 39916800)^{tr}, (1060, 97308, 1494576, 3265920, 0)^{tr},\\ (620, 26908, 84816, -282240, 0)^{tr}, (290, 1828, -30024, 50400, 0)^{tr},\\ (70, -1692, 9576, -12960, 0)^{tr}, (-40, 508, -2304, 2880, 0)^{tr}}$.
\item The degree of $f$ is 12 and the following inequalities are satisfied.
\begin{itemize}
\item[(i)] $\scalar{(0, 0, 0, 0, 0, 1),\mathbf{v}^{12}}\geq 0$
\item[(ii)] $\scalar{(16, 8, 4, 2, 1, -6),\mathbf{v}^{12}}+32\geq 0$
\item[(iii)] $\scalar{(3680, 544, 56, -8, -10, 33),\mathbf{v}^{12}}+22912\geq 0$
\item[(iv)] $\scalar{(75120, 2568, -228, -6, 45, -110),\mathbf{v}^{12}}+1192224\geq 0$
\item[(v)] $\scalar{(187520, -4352, 32, 64, -80, 165),\mathbf{v}^{12}}+9601024\geq 0$
\item[(vi)] $\scalar{(-7680, -984, 144, -42, 35, -66),\mathbf{v}^{12}}+9717408\geq 0$
\item[(vii)] $\scalar{(-104672, 3104, -248, 56, -42, 77),\mathbf{v}^{12}}+15276416\geq 0$
\end{itemize}
In particular, $\mathbf{v}^{12}$ lies in the convex hull of

$\set{(2306, 530368, 25971144, 302788656, 723263040, 239500800)^{tr}, \\(1580, 232708, 6523056, 34999200, 21772800, 0)^{tr},\\ (986, 78268, 804024, -840384, -2177280, 0)^{tr}, (524, 13588, -95952, -108576, 483840, 0)^{tr}, \\(194, -2912, -2232, 105264, -181440, 0)^{tr}, (-4, -932, 15984, -80064, 103680, 0)^{tr}, \\(-70, 1708, -17544, 72000, -86400, 0)^{tr}}$.
\item The degree of $f$ is 13 and the following inequalities are satisfied.
\begin{itemize}
\item[(i)] $\scalar{(0, 0, 0, 0, 0, 1),\mathbf{v}^{13}}\geq 0$
\item[(ii)] $\scalar{(96, 48, 24, 12, 6, -11),\mathbf{v}^{13}}+192\geq 0$
\item[(iii)] $\scalar{(6256, 968, 124, 2, -9, 9),\mathbf{v}^{13}}+38432\geq 0$
\item[(iv)] $\scalar{(86160, 4200, -60, -30, 15, -11),\mathbf{v}^{13}}+1260960\geq 0$
\item[(v)] $\scalar{(1583520, -3984, -1560, 300, -90, 55),\mathbf{v}^{13}}+60906432\geq 0$
\item[(vi)] $\scalar{(2849920, -117280, 6880, -760, 180, -99),\mathbf{v}^{13}}+590783360\geq 0$
\item[(vii)] $\scalar{(-1475936, 25592, -1064, 98, -21, 11),\mathbf{v}^{13}}+344354528\geq 0$
\end{itemize}
In particular, $\mathbf{v}^{13}$ lies in the convex hull of

$\set{(3206, 1071868, 81406344, 1627117056, 7827719040, 6227020800)^{tr},\\ (2270, 508708, 23753736, 249815520, 439084800, 0)^{tr}, \\(1490, 196708, 4275576, 6979680, -32659200, 0)^{tr}, \\(866, 50068, -22536, -2702304, 5080320, 0)^{tr}, \\(398, 148, -123000, 847008, -1209600, 0)^{tr}, \\(86, -4532, 61704, -287424, 362880, 0)^{tr}, \\(-70, 1708, -17544, 72000, -86400, 0)^{tr}}$.
\item The degree of $f$ is 14 and the following inequalities are satisfied.
\begin{itemize}
\item[(i)] $\scalar{(0, 0, 0, 0, 0, 0, 1),\mathbf{v}^{14}}\geq 0$
\item[(ii)] $\scalar{(32, 16, 8, 4, 2, 1, -7),\mathbf{v}^{14}}+64\geq 0$
\item[(iii)] $\scalar{(45696, 7296, 1056, 96, -24, -24, 91),\mathbf{v}^{14}}+278016\geq 0$
\item[(iv)] $\scalar{(573216, 33888, 744, -168, 6, 33, -91),\mathbf{v}^{14}}+7911552\geq 0$
\item[(v)] $\scalar{(24080000, 269440, -22240, 1120, 200, -440, 1001),\mathbf{v}^{14}}+771857920\geq 0$
\item[(vi)] $\scalar{(30693600, -533520, 9720, 540, -450, 495, -1001),\mathbf{v}^{14}}+3152571840\geq 0$
\item[(vii)] $\scalar{(-7032576, -169728, 15936, -2112, 624, -528, 1001),\mathbf{v}^{14}}+6989196288\geq 0$
\item[(viii)] $\scalar{(-5559008, 96992, -4088, 392, -98, 77, -143),\mathbf{v}^{14}}+1292719232\geq 0$
\end{itemize}
In particular, $\mathbf{v}^{14}$ lies in the convex hull of

$\set{(4347, 2046870, 231190132, 7475341368, 67360465152, 138619313280, 43589145600)^{tr},\\ (3164, 1036588, 76136688, 1431621504, 6044647680, 3353011200, 0)^{tr}, \\(2163, 445998, 17610220, 114649848, -163634400, -279417600, 0)^{tr}, \\(1344, 142968, 1436608, -15125616, -2062080, 50803200, 0)^{tr}, \\(707, 19390, -492228, 1431288, 6560352, -15240960, 0)^{tr}, \\(252, -8820, 51952, 524928, -4578048, 6773760, 0)^{tr}, \\(-21, -1722, 66148, -787656, 3441600, -4233600, 0)^{tr}, \\(-112, 4648, -89280, 808848, -3110400, 3628800, 0)^{tr}}$.
\end{itemize}

\end{document}